\documentclass[final,5p,times,twocolumn]{elsarticle}
\pdfoutput=1

\usepackage{amsmath}
\usepackage{amssymb}
\usepackage{amsthm}
\usepackage{graphicx}
\usepackage{xcolor}
\usepackage{microtype}
\usepackage{url}

\DeclareMathOperator{\tr}{tr}

\theoremstyle{plain}
\newtheorem{maintheorem}{Theorem}
\newtheorem{proposition}{Proposition}
\theoremstyle{definition}
\newtheorem{definition}{Definition}
\newtheorem{remark}{Remark}

\journal{Statistics \& Probability Letters}

\begin{document}
\begin{frontmatter}

\title{Largest Finite Root of Identity-Scale Doubly Singular Beta Type II Ensembles}

\author[inst1]{Stepan Grinek\corref{cor1}}
\ead{stepan.grinek@jax.org}
\cortext[cor1]{Corresponding author}
\address[inst1]{The Jackson Laboratory-New York Stem Cell Foundation Research Institute, New York, NY, USA.}

\begin{abstract}
Classical largest-root distributions for Wishart ratios and matrix-variate
beta ensembles are usually formulated when the denominator Wishart matrix
is nonsingular. In many high-dimensional settings, however, the ambient dimension $p$ exceeds both Wishart degrees of freedom, so the corresponding beta ensemble is doubly singular and the usual matrix product $A^{-1}B$ is not defined. We consider independent central Wishart matrices $A\sim W_p(m,I_p)$ and $B\sim W_p(q,I_p)$ in the regime $p>m\ge q$. For the finite generalized roots of the pair $(B,A)$ or, equivalently, the
nonzero eigenvalues of $BA^+$, we prove the exact identity
\begin{equation*}
\lambda_{\max}
\stackrel{d}{=}
\lambda_{\max}\left\{W_q(m,I_q)W_q(p-m+q,I_q)^{-1}\right\}.
\end{equation*}
Here $W_d(r,I_d)$ denotes a $d\times d$ central Wishart matrix with
$r$ degrees of freedom and identity scale. Thus, the identity-scale
doubly singular beta type II largest-root problem reduces exactly to a
nonsingular $q$-dimensional Roy statistic, making its finite-sample CDF
directly accessible to classical largest-root algorithms.
\end{abstract}


\end{frontmatter}

\section{Introduction}

Generalized roots of pairs of independent Wishart matrices are central to
MANOVA, Roy's largest-root test, and the theory of matrix-variate beta and
Jacobi ensembles. In the classical nonsingular beta type II setting, one
studies the eigenvalues of $A^{-1}B$, where $A$ and $B$ are independent
Wishart matrices with a common scale matrix. The largest root has a long history in multivariate statistics and random
matrix theory; see, for example, Mardia et al.~\cite{Mardia} and
Johnstone~\cite{Johnstone}. 

A different regime occurs when the ambient dimension $p$ exceeds both Wishart degrees of freedom. Then $A$ and $B$ are singular and
the usual inverse $A^{-1}$ does not exist. If both matrices have rank
smaller than $p$, the associated beta ensemble is doubly singular.
When $A$ is singular, we define the finite roots by restricting the
generalized eigenvalue problem to $\mathcal R(A)$, the range of $A$.
Section~\ref{sec:notation} shows that these are the positive eigenvalues of
the compressed matrix in Definition~\ref{def:finite-roots}, equivalently
the nonzero eigenvalues of $BA^+$.

Doubly singular beta type I and type II distributions, including the
joint densities of their nonzero eigenvalues, were developed by
D\'iaz-Garc\'ia and Guti\'errez-J\'aimez \cite{Diaz,Diaz2}. Related
large-$p$ approximations were studied by Srivastava and Fujikoshi
\cite{sri1,sri2,sri3}, and Tracy--Widom limits for Jacobi largest roots
were obtained by Johnstone \cite{Johnstone}. Exact CDF representations
for singular beta-Wishart and beta $F$ problems have subsequently been
derived using heterogeneous hypergeometric functions with two matrix
arguments \cite{ShimizuHashiguchi2021,ShimizuHashiguchi2022}.

The present contribution is an exact distributional reduction in the
central identity-scale case: the largest finite root is identified with
a lower-dimensional, nonsingular Roy statistic. This avoids direct evaluation of hypergeometric functions with two
matrix arguments and makes classical
double-Wishart algorithms immediately applicable.

The proof uses the orthogonal invariance of standard Gaussian matrices to
remove the random orientation of $\mathcal R(A)$, followed by the classical
dimension-duality identity for Roy's largest root. The resulting
representation is finite-sample exact and places the doubly singular problem
within the standard Roy framework.

An applied motivation comes from high-dimensional multivariate analysis, including Principal Components of Heritability (PCH) and related MANOVA-type statistics in genomic studies \cite{PCH}. In such applications, the number of measured variables can exceed the number of observations, making both model and residual covariance matrices singular. The result below supplies an exact identity-scale reference distribution for the largest finite root and a benchmark for large-$p$ approximations.

The R code used to generate Figure~\ref{fig:numerics} is available at
\url{https://github.com/stepanv1/DoubleWishart/blob/master/DW.R}.

\section{Notation and finite roots}\label{sec:notation}

We write $W_d(r,\Sigma)$ for a $d\times d$ central Wishart matrix with
$r$ degrees of freedom and scale matrix $\Sigma$. The notation
$\stackrel{d}{=}$ denotes equality in distribution, while
$\xrightarrow{d}$ and $\xrightarrow{\mathbb P}$ denote convergence in
distribution and in probability, respectively. For the matrices considered below, whose relevant eigenvalues are real,
$\lambda_{\max}(M)$ denotes the largest eigenvalue of $M$.

Let $A\sim W_p(m,\Sigma)$ and $B\sim W_p(q,\Sigma)$ be independent
central Wishart matrices with $p>m\ge q$ and common scale matrix $\Sigma$.
Write the rank-$m$ spectral decomposition of $A$ as
\begin{equation*}
A=HLH^T,
\end{equation*}
where $H\in\mathbb R^{p\times m}$ satisfies $H^TH=I_m$ and $L$ is
diagonal with positive entries. Define the whitening map
\begin{equation*}
T=HL^{-1/2}.
\end{equation*}
Then $T^TAT=I_m$.

We now whiten $A$ on its range and solve the resulting ordinary
eigenvalue problem; the equivalence with the projected generalized equation
is verified immediately below.

\begin{definition}[Finite roots of the doubly singular beta type II ensemble]
\label{def:finite-roots}
The finite roots of the pair $(B,A)$ are the positive eigenvalues of
\begin{equation}
K=T^TBT=L^{-1/2}H^TBHL^{-1/2}.
\label{eq:finite-root-matrix}
\end{equation}
Their largest value is denoted by $\lambda_{\max}$.
\end{definition}

The connection with the generalized eigenvalue problem requires some
care because $A$ is singular. We restrict the problem to
$\mathcal R(A)$, the range of $A$, and project the equation back onto
that subspace. Every $e\in\mathcal R(A)$ can be written uniquely as
$e=Hv$, with $v\in\mathbb R^m$. The projected equation is
\begin{equation}
H^T(Be-\lambda Ae)=0,
\qquad e=Hv,
\label{eq:projected-generalized-problem}
\end{equation}
or equivalently
\begin{equation*}
H^TBH\,v=\lambda Lv.
\end{equation*}
Setting $u=L^{1/2}v$, or equivalently, $e=Tu$, transforms this equation
into
\begin{equation*}
Ku=\lambda u.
\end{equation*}
Thus, the eigenvalues of $K$ are precisely the generalized roots of the
pencil after restriction to $\mathcal R(A)$. Equivalently, they are the
stationary values of
\begin{equation*}
\frac{e^TBe}{e^TAe},
\qquad e\in\mathcal R(A)\setminus\{0\}.
\end{equation*}

The construction in \eqref{eq:finite-root-matrix} is the compressed
form of the Moore--Penrose definition used for singular and doubly
singular beta type II matrices \cite{Diaz, Diaz2}. Indeed,
\begin{equation*}
(A^+)^{1/2}=(A^{1/2})^+=HL^{-1/2}H^T
\end{equation*}
and hence
\begin{equation}
(A^+)^{1/2}B(A^+)^{1/2}=HKH^T.
\label{eq:embedded-beta-II}
\end{equation}
Because $H^TH=I_m$, the matrices in
\eqref{eq:finite-root-matrix} and \eqref{eq:embedded-beta-II} have the
same positive eigenvalues, including multiplicities. 

The same roots can also be represented as the nonzero eigenvalues of
$BA^+$. Since
\begin{equation*}
A^+=HL^{-1}H^T,
\end{equation*}
define
\begin{equation*}
P=BHL^{-1/2},
\qquad
Q=L^{-1/2}H^T.
\end{equation*}
Then
\begin{equation*}
PQ=BA^+,
\qquad
QP=K.
\end{equation*}
The products $PQ$ and $QP$ have the same nonzero eigenvalues, including
algebraic multiplicities. Therefore, the positive eigenvalues of $K$
coincide exactly with the nonzero eigenvalues of $BA^+$. Although
$BA^+$ is generally nonsymmetric, these eigenvalues are real and
nonnegative because $K$ is symmetric positive semidefinite.

To connect the compressed matrix $K$ with the definitions used in the singular-beta literature, write $B=YY^T$, where $Y\in\mathbb R^{p\times q}$. One of the Moore--Penrose formulations of the doubly singular beta type II
matrix considered by D\'iaz-Garc\'ia and Guti\'errez-J\'aimez
\cite{Diaz,Diaz2} is
\begin{equation*}
\widetilde F=Y^TA^+Y.
\end{equation*}
With
\[
G=L^{-1/2}H^TY\in\mathbb R^{m\times q},
\]
we have
\[
K=GG^T,
\qquad
\widetilde F=G^TG.
\]
The matrices $GG^T$ and $G^TG$ have the same positive eigenvalues,
including multiplicities. Consequently, the finite roots defined through $K$ coincide samplewise, and hence also in distribution, with the nonzero roots in the corresponding Moore--Penrose definition of the doubly singular beta type II ensemble. Since $\operatorname{rank}(X)\le q$, there are at most $q$ positive roots; if $q<m$, the remaining eigenvalues of the $m\times m$ matrix $K$ are zero.

\section{Exact identity-scale representation}\label{sec:exact}

We first state the exact finite-dimensional distributional result, which is the central contribution. Large-$p$ approximations and scale corrections are discussed afterwards in Section~\ref{sec:asymptotic}.
\begin{maintheorem}[Exact double-Wishart representation]
\label{thm:exact}
Let $A\sim W_p(m,I_p)$ and $B\sim W_p(q,I_p)$ be independent central Wishart matrices with $p>m\ge q$. Let $\lambda_{\max}$ denote the largest finite generalized root of the pair $(B,A)$ in the sense of Definition~\ref{def:finite-roots}. Then
\begin{equation}\label{eq:m-representation}
\lambda_{\max}
\stackrel{d}{=}
\lambda_{\max}\left\{W_m(q,I_m)W_m(p,I_m)^{-1}\right\},
\end{equation}
where the right-hand side is understood as the largest nonzero root of the corresponding double-Wishart product. Equivalently,
\begin{equation}\label{eq:q-representation}
\lambda_{\max}
\stackrel{d}{=}
\lambda_{\max}\left\{W_q(m,I_q)W_q(p-m+q,I_q)^{-1}\right\}.
\end{equation}
The Wishart matrices appearing on the right-hand sides are independent.
\end{maintheorem}

\begin{proof}
Use independent standard Gaussian matrices
$X\in\mathbb R^{p\times m}$ and $Z\in\mathbb R^{p\times q}$ to write
$A=XX^T$ and $B=ZZ^T$. Write the singular-value decomposition using the $m$ nonzero singular
values of $X$ as
\[
  X=HL^{1/2}R^T,
\]
where $H^TH=I_m$, $R^TR=I_m$, and $L$ contains the positive
eigenvalues of $A$.  Then
$A=HLH^T$ and $X^TX=RLR^T\sim W_m(p,I_m)$.

By Definition~\ref{def:finite-roots}, the finite roots are the positive
eigenvalues of
\[
  L^{-1/2}H^TZZ^THL^{-1/2}.
\]
If $V=H^TZ$, the nonzero eigenvalues of this matrix are those of
\[
  V^TL^{-1}V.
\]
Conditional on $X$, the matrix $H$ is fixed and the columns of
$V=H^TZ$ are independent $N_m(0,I_m)$ vectors, because $H^TH=I_m$.
This conditional distribution is the same for every value of $X$.
Therefore, $V$ is independent of $X$, and hence independent of $L$ and $R$.

Set
\[
  S=X^TX=RLR^T,
  \qquad W=RV.
\]
Conditional on $X$, the matrix $R$ is fixed and orthogonal, so
$W=RV$ is again standard Gaussian. Its conditional law does not depend on
$X$; consequently, $W$ is independent of $X$ and hence of $S=X^TX$.  Moreover,
\[
  W^TS^{-1}W
  =V^TR^T(RL^{-1}R^T)RV
  =V^TL^{-1}V.
\]
The nonzero eigenvalues of $W^TS^{-1}W$ are the same as those of
$WW^TS^{-1}$.  Since
\[
  WW^T\sim W_m(q,I_m),
  \qquad S\sim W_m(p,I_m),
\]
and these matrices are independent, this proves
\eqref{eq:m-representation}.  The product $WW^TS^{-1}$ is similar to the
symmetric positive semidefinite matrix
$S^{-1/2}WW^TS^{-1/2}$, so its nonzero roots are real and nonnegative.

Finally, the standard dimension-duality identity for Roy's largest root
\cite{Mardia,Chiani} gives
\begin{align*}
&\lambda_{\max}\left\{
  W_m(q,I_m)W_m(p,I_m)^{-1}
 \right\}\\
&\qquad\stackrel d=
 \lambda_{\max}\left\{
  W_q(m,I_q)W_q(p-m+q,I_q)^{-1}
 \right\}.
\end{align*}
This proves \eqref{eq:q-representation}.
\end{proof}

\begin{remark}[Beta type I scale]
Some algorithms and tables for Roy's largest root are written for beta type I, or Jacobi, roots $\theta\in(0,1)$ rather than beta type II roots $\lambda\in(0,\infty)$. The two parameterizations are related by
\begin{equation*}
\theta=\frac{\lambda}{1+\lambda},
\qquad
\lambda=\frac{\theta}{1-\theta}.
\end{equation*}
Thus, a CDF evaluated on the Jacobi scale can be transferred to the beta type II scale by this monotone transformation.
\end{remark}

The representation \eqref{eq:q-representation} is the computationally useful form of the theorem. It reduces the original $p$-dimensional doubly singular problem to a $q$-dimensional nonsingular double-Wishart largest-root problem. The CDF can then be computed using available algorithms for the classical largest root distribution, such as the method of Chiani \cite{Chiani}.
\section{Large-p approximations and scale effects}\label{sec:asymptotic}

The exact reduction exploits the spherical symmetry of the
identity-scale Gaussian model. For $\Sigma\ne I_p$, the orientation of
$\mathcal R(A)$ interacts with $\Sigma$, so the identity-scale argument no
longer applies. A simultaneous orthogonal change of coordinates preserves the
generalized roots while replacing $\Sigma$ by $U^T\Sigma U$; hence, their
distribution depends on $\Sigma$ only through its eigenvalues. The following
large-$p$ approximations describe the leading effect of a common nonidentity
scale and connect the exact result with earlier singular MANOVA theory.

\begin{proposition}[Identity-scale large-$p$ approximation]\label{prop:identity-asymptotic}
Let $A\sim W_p(m,I_p)$ and $B\sim W_p(q,I_p)$ be independent central Wishart matrices with $m$ and $q$ fixed while $p\to\infty$. Let $\lambda_{\max}$ be the largest finite generalized root of $(B,A)$. Then
\begin{equation}\label{eq:identity-asymptotic}
p\lambda_{\max}
\xrightarrow{d}
\lambda_{\max}\{W_q(m,I_q)\}.
\end{equation}
\end{proposition}
This follows immediately from Theorem~\ref{thm:exact}, since, for fixed
$m$ and $q$,
\begin{equation*}
p^{-1}W_q(p-m+q,I_q)\xrightarrow{\mathbb P}I_q.
\end{equation*}
The continuous mapping theorem and Slutsky's theorem then give
\eqref{eq:identity-asymptotic}. This recovers the identity-scale
large-$p$ approximation used in the singular MANOVA literature
\cite{sri1,sri2}.

\paragraph{Common nonidentity scale}
Let $\Sigma=\Sigma_p$ and suppose that $p\to\infty$ with $m$ and $q$ fixed.
Assume the conditions of Srivastava and Fujikoshi
\cite[Eq.~(2.7)]{sri2}; in particular, let
\[
a_{i,p}=\frac{\tr(\Sigma_p^i)}{p}\longrightarrow a_i\in(0,\infty),
  \qquad i=1,2,
\]
and define
\[
  \beta_\Sigma=\frac{a_1^2}{a_2}.
\]
Write $A=HLH^T$ and $B=ZZ^T$, where the columns of $Z$ are
independent $N_p(0,\Sigma_p)$ vectors, and set
$C=(H^T\Sigma_p H)^{-1/2}$ and $V=CH^TZ$. Then the nonzero
eigenvalues of $BA^+$ coincide with the eigenvalues of
$V^T(CLC)^{-1}V$. Conditional on $A$, the columns of $V$ are
independent $N_m(0,I_m)$ vectors. The
projected-covariance limit of Srivastava and Fujikoshi
\cite[Eq.~(2.7)]{sri2} gives
\[
  p(CLC)^{-1}\xrightarrow{\mathbb P}\frac{1}{\beta_\Sigma}I_m.
\]
Because the conditional law of $V$ does not depend on $A$, $V$ is
independent of $CLC$. Slutsky's theorem therefore yields
\[
  p\lambda_{\max}
  \xrightarrow{d}
  \lambda_{\max}\left\{
    \frac{1}{\beta_\Sigma}W_q(m,I_q)
  \right\}.
\]

In the separate regime $m,p\to\infty$, Srivastava~\cite{sri3}
proposed the consistent estimators
\begin{align*}
\hat a_1 &= \frac{\tr A}{mp},\\
\hat a_2 &=
\frac{\tr(A^2)-m^{-1}\{\tr A\}^2}
     {(m-1)(m+2)p}.
\end{align*}
The corresponding plug-in estimator of the scale factor is
$\widehat\beta_\Sigma=\hat a_1^2/\hat a_2$.
\section{Numerical illustration}\label{sec:numerics}
We compare the empirical CDF from the original $p$-dimensional
identity-scale model with the exact $q$-dimensional double-Wishart CDF in
\eqref{eq:q-representation} and Johnstone's Jacobi-ensemble Tracy--Widom
approximation \cite{Johnstone}.
Each replicate draws independent standard Gaussian matrices
$X\in\mathbb R^{p\times m}$ and $Y\in\mathbb R^{p\times q}$, with
$A=XX^T$ and $B=YY^T$, exactly matching the model of
Theorem~\ref{thm:exact}.

CDFs are evaluated for the beta type II root $\lambda_{\max}$ and plotted
against $\log_{10}(\lambda_{\max})$. For the exact
calculation, the root is transformed to the Jacobi scale
\[
  \theta=\frac{\lambda}{1+\lambda},
\]
and the canonical Roy parameters are
\[
  s=q,\qquad
  m_{\mathrm C}=\frac{m-q-1}{2},\qquad
  n_{\mathrm C}=\frac{p-m-1}{2}.
\]
The exact CDF is evaluated with \texttt{rootWishartHD}, using
\texttt{doubleWishart} in regular cases and the arbitrary-precision
\texttt{doubleWishart\_robustPfaffians} backend near the boundary or for large
$q$; the stress case uses the latter explicitly.
The Tracy--Widom curve uses Johnstone's logit-Jacobi centering and scaling,
followed by the real Tracy--Widom CDF \texttt{RMTstat::ptw} with $\beta=1$.

The primary calculation fixes $m=100$, $q=6$, and uses
$p\in\{101,110,120,130,200,500\}$.
A separate stress case uses $(p,m,q)=(100,99,98)$. Under the exact reduction,
this case becomes $G_1G_2^{-1}$, where $G_1$ and $G_2$ are independent
and
\[
G_i\sim W_{98}(99,I_{98}),\qquad i=1,2.
\]
Thus both factors have only one degree of freedom beyond their dimension. The stress case therefore lies
near the boundary of the nonsingular Roy parameter space and provides a
stringent comparison with the Tracy--Widom approximation.

\begin{figure*}[t]
\centering
\includegraphics[width=\textwidth]{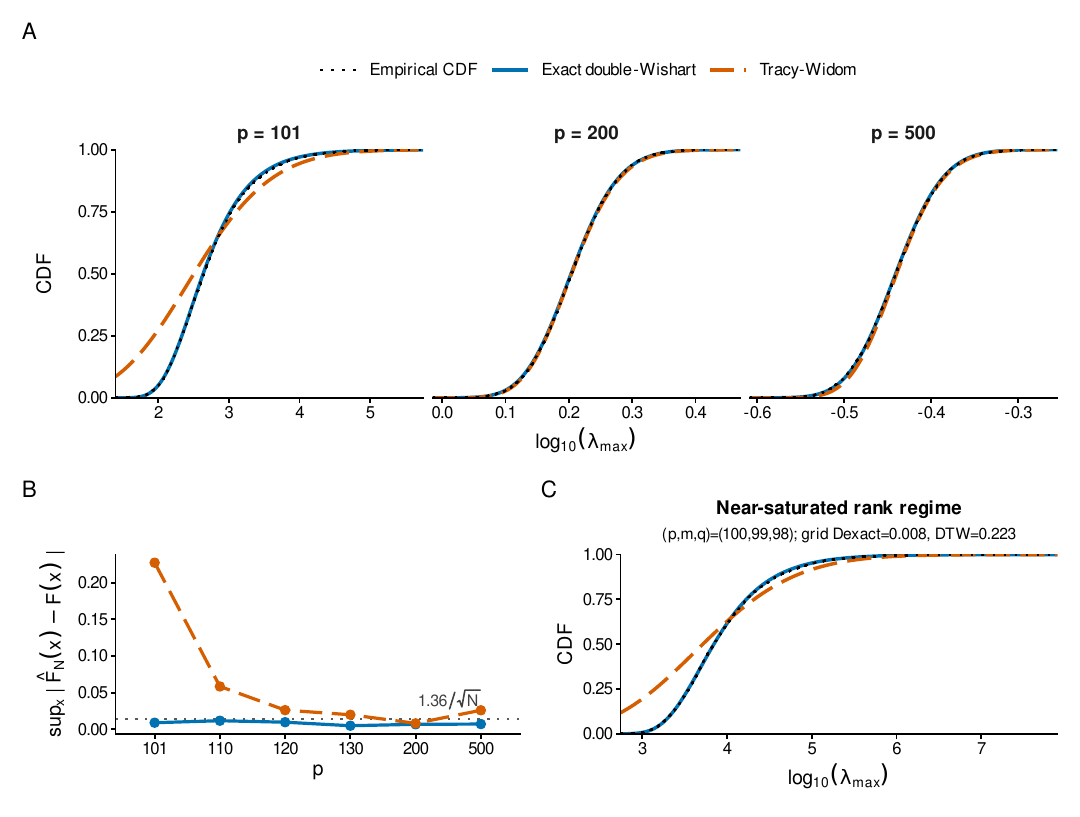}
\caption{\textbf{Finite-sample validation and asymptotic comparison.}
\textbf{(A)} Empirical CDFs (black dotted) of the largest finite root from the
original $p$-dimensional Gaussian construction, together with the exact
reduced double-Wishart CDF (blue solid) and the Tracy--Widom approximation
(orange dashed). Results use $m=100$, $q=6$, identity scale,
$p\in\{101,200,500\}$, and $N=10{,}000$ Monte Carlo replicates for each value
of $p$; the horizontal axis is $\log_{10}(\lambda_{\max})$.
\textbf{(B)} Maximum absolute CDF discrepancy over
$p\in\{101,110,120,130,200,500\}$. The dotted horizontal line is the nominal
5\% one-sample Kolmogorov--Smirnov reference
$1.36/\sqrt{N}=0.0136$.
\textbf{(C)} Near-saturated case $(p,m,q)=(100,99,98)$ based on
$N_{\mathrm{stress}}=5{,}000$ replicates. The grid discrepancies are
$D_{\mathrm{exact}}=0.0083$ and $D_{\mathrm{TW}}=0.2228$, compared with the
Monte Carlo reference
$1.36/\sqrt{N_{\mathrm{stress}}}=0.0192$.}
\label{fig:numerics}
\end{figure*}

For a candidate CDF $F$, define
\begin{equation*}
  D(F)=\sup_x\left|\widehat F_N(x)-F(x)\right|.
\end{equation*}
For a continuous candidate CDF, the one-sample Kolmogorov--Smirnov
distance is attained immediately before or at an observed root; both sides of
each empirical jump are therefore checked in the primary calculation.
Figure~\ref{fig:numerics}A shows representative CDFs, while
Figure~\ref{fig:numerics}B reports
$D_{\mathrm{exact}}$ and $D_{\mathrm{TW}}$ on all dimensions tested. The
stress-case discrepancies are evaluated on the common 801-point plotting
grid. The line $1.36/\sqrt{N}$ is shown as the usual asymptotic 5\% reference
for the Monte Carlo sampling variation.

The exact reduced distribution agrees with simulation at Monte Carlo
accuracy. Across the six primary values of $p$, $D_{\mathrm{exact}}$ ranges
from $0.0048$ to $0.0116$, below the reference value $0.0136$. The Tracy--Widom discrepancy depends strongly on the excess degrees of
freedom of the two reduced Wishart factors. It is $0.2274$ when
$p-m=1$ at $p=101$, decreases to $0.0584$, $0.0260$, and $0.0197$ at
$p=110$, $120$, and $130$, respectively, and equals $0.0080$ and $0.0258$
at $p=200$ and $p=500$. In the near-saturated case, where both
$p-m=1$ and $m-q=1$, the exact grid discrepancy is $0.0083$, below the
Monte Carlo reference $0.0192$, whereas the Tracy--Widom discrepancy is
$0.2228$. These results confirm the finite-sample reduction in both regular
and boundary configurations and show that a large ambient dimension alone
does not ensure an accurate Tracy--Widom approximation.

\section{Discussion}\label{sec:discussion}

Theorem~\ref{thm:exact} identifies the largest finite root of the
identity-scale doubly singular beta type II ensemble with the largest
root of a nonsingular $q$-dimensional Roy statistic. The reduction
converts the original singular $p$-dimensional problem into a classical
double-Wishart problem in dimension $q$, making the finite-sample CDF
directly accessible to established largest-root algorithms. Thus it provides both an exact distributional characterization and a
practical computational route.

Figure~\ref{fig:numerics} shows that the exact reduced CDF tracks the
simulated distribution in both regular and near-saturated regimes. By
contrast, the Tracy--Widom approximation can deteriorate sharply near the
boundary of the classical Roy parameter space. The exact reduction is especially useful in configurations where the
asymptotic approximation is least reliable.

The exact identity uses the spherical symmetry of the identity-scale
model. For a common nonidentity scale, Section~\ref{sec:asymptotic} provides
the corresponding large-$p$ approximation; at finite $p$, the distribution
still depends on the alignment between the random range $\mathcal R(A)$
and $\Sigma$.

\end{document}